\documentclass[letterpaper, 10 pt, conference]{ieeeconf}  

\IEEEoverridecommandlockouts                              

\title{\LARGE \bf
Indirect controllability of two interacting qubits in presence of dissipation: a first analysis.
}

\author{K. Verzhanska$^{1}$ and F. C. Chittaro$^{1}$
\thanks{*This work was supported by QUACO, PRC ANR-17-CE40-0007-01 and by CARTT-IUT de Toulon}
\thanks{$^{1}$LIS, UMR CNRS 7020, Universit\'e de Toulon, Aix Marseille Univ., France
        {\tt\small  kateryna.verzhanska@lis-lab.fr francesca.chittaro@univ-tln.fr}}%
}

\usepackage[utf8]{inputenc}
\usepackage{amsmath,amsfonts,amssymb}
\usepackage{enumerate}
\usepackage{mathtools}
\usepackage{xcolor}
\mathtoolsset{showonlyrefs,showmanualtags}

\usepackage{graphicx}
\newcommand*{\Scale}[2][4]{\scalebox{#1}{$#2$}}%

\newcommand{\bq}{\boldsymbol{q}}
\newcommand{\bu}{\boldsymbol{u}}

\newcommand{\vv}{\boldsymbol{v}}

\newcommand{\Tr}{\mathrm{Tr}}
\newcommand{\id}{\boldsymbol{I}}

\newcommand{\llangle}{\langle\!\langle}
\newcommand{\rrangle}{\rangle\!\rangle}
\newcommand{\brho}{\boldsymbol{\rho}}

\newcommand{\Hdue}[1]{h_{#1}}
\newcommand{\Htre}[1]{\hat{h}_{#1}}
\newcommand{\Hse}[1]{M_{#1}}
\newcommand{\Hqu}[1]{\widehat{H}_{#1}}
\newcommand{\bvz}{\mathbf{v}_0}

\newtheorem{remark}{Remark}
\newtheorem{proposition}{Proposition}

\newtheorem{theorem}{Theorem}
\newtheorem{corollary}{Corollary}
\newtheorem{prob}{Problem}

\begin{document}

\maketitle
\thispagestyle{empty}
\pagestyle{empty}

\begin{abstract}
We consider a bipartite open quantum system constituted by two interacting qubits $A$ and $B$, assuming that the former is coupled to the environment and is directly affected by coherent control, while the latter does not interact directly with the environment and the control fields.
We are interested in the controllability properties of the subsystem $B$.

In this paper, we give a first analysis of the problem and provide some negative answers.
\end{abstract}

\section{Introduction}

Quantum control deals with the 
manipulation of dynamical systems at the molecular and atomic scale, where the dynamics are governed by quantum mechanics.
As in classical control, the notion of controllability refers to the possibility to steer a given initial state to any desired target
state, by applying appropriate external fields. 

For finite-dimensional closed quantum systems, the Lie Algebraic Rank Condition (LARC) is a necessary and sufficient conditions for controllability of the bilinear Schrödinger equation (\cite{albertini-dalessandro,altafiniJMP}). For open quantum systems, the problem is more delicate: in \cite{altafiniJMP}, the author proved that an open quantum system is never small-time locally controllable (STLC) by means of coherent control, and that some configurations are not reachable in finite time; in \cite{dirr-helmke}, it has been pointed out that, as the underlying group is not compact, then LARC is a sufficient condition for accessibility, but does not imply controllability. 

Yet, in many experimental set-ups, there is no need to control the whole system: 
for instance, in typical situations the system of interest ($B$) is well isolated from the environment, and interacts with an (eventually open) {\it accessory system} (or {\it ancilla}), which can be directly controlled  (\cite{aiello-cappellaro15}, \cite{burgarth2010scalable}, \cite{dyn-cat}). In such situations, a natural question is to analyze the controllability properties of the subsystem $B$ only; this notion is called {\it indirect controllability}.

For closed quantum systems, a detailed analysis of indirect controllability has been carried out in the papers~\cite{dal-rom,alb-rom-dal}, in terms of Lie groups theory. 
For open quantum systems, the problem has been investigated in specific situations (see for instance~\cite{grimmer,layden}) but, to our knowledge, a general theoretical analysis is missing. 

In this paper, we focus on a composite quantum system made by two interacting qubits $A$ and $B$, such that
\begin{itemize}
\renewcommand{\labelitemi}{$\cdot$}
\item the control acts directly only on the subsystem $A$, which is subject to the interaction with the environment;
\item the system of interest $B$ does not interact directly with the controls and the environment (only through its interaction with $A$).
\end{itemize}

As a first step towards the characterization of indirect controllability for such systems, we restrict our attention to a particular class of target states: the states (of the whole system) such that their reduction to $B$ is ``pure'' (\cite{nielsen-chuang}). In particular, we ask ourselves the following question: is it possible to ``purify''  the state of the subsystem $B$ and/or to keep it pure? 

Even if the existence of purifying dynamics is well known (\cite{altafiniJMP}), they achieve complete purification only asymptotically in time. Based on this fact, we remark in Corollary~\ref{cor: ppure infinite time} that also ``partial purification'' is only asymptotic.

We then investigate the possibility of protecting the subsystem $B$ from dissipation, that is, to keep its state pure. We analyze three possible interactions between the qubits, and show that the only dynamics that conserve partial purity are trivial or not affected by the control. 

These partial results provide a first negative answer to the question of indirect controllability in presence of dissipation.

\medskip
The structure of the paper is the following: 
in Sec.~\ref{sec: preliminaries} we introduce the minimal relevant notions on bipartite quantum systems, and we exhibit the class of systems we are interested into; in Sec.~\ref{sec: K}, we discuss the structure of the set of admissible (physical) states, with a particular interest on states that correspond to pure reduced states; finally, in Section~\ref{sec: results} we derived some
preliminary results.


\section{Statement of the problem} \label{sec: preliminaries}
\subsection{Basic facts on open quantum systems}
In this section, we just provide a minimal description of the formalism of open quantum systems. For more details, we refer to the monographs~\cite{breuer-petr,haroche,nielsen-chuang}.

Let us first start with some notations:
we denote with $\mathfrak{her}(n)$ (respectively, $\mathfrak{her}_0(n)$)  the set of Hermitian (respectively, traceless) $n$ dimensional  matrices, with 
$\mathfrak{su}(n)$ the set of anti-Hermitian traceless $n$ dimensional  matrices 
and with $\mathfrak{so}(n)$  the set of anti-symmetric real $n$ dimensional matrices.

Consider the Hilbert space $\mathbb{C}^{n_A}\otimes \mathbb{C}^{n_B}$, for $n_A,n_B\geq 2$. The {\it partial trace over $\mathbb{C}^{n_A}$} is the unique linear operator  
$\Tr_A: \mathfrak{her}(n_An_B) \to \mathfrak{her}(n_B)$ 
such that for every $M_A\in\mathfrak{her}(n_A)$ and every $M_B\in\mathfrak{her}(n_B)$ it holds $\Tr_A (M_A\otimes M_B)=M_B(\Tr M_A)$, where $\Tr$ denotes the usual trace operation on $\mathfrak{her}(n_A)$. The partial trace over $\mathbb{C}^{n_B}$ is defined analogously.

In the standard formulation of quantum mechanics, the state of a finite-dimensional quantum system is represented by a positive semi-definite Hermitian operator with trace one, acting on a finite-dimensional complex Hilbert space $\mathcal{H}$ (usually identified with $\mathbb{C}^N$, for some $N$);   
such operator is called {\it density matrix} (or {\it density operator}) and is usually denoted by $\brho$. 
In this paper, we denote with $\mathcal{P}$ the set of density matrices (i.e., positive semi-definite and Hermitian with trace one) on $\mathcal{H}$; 
$\mathcal{P}$ is a compact and convex subset of $\mathfrak{her}(N)$, and its extreme points coincide with rank one projection operators on $\mathbb{C}^N$; these states, called {\it pure states} in the language of quantum mechanics, are characterized by the the property $\Tr \brho^2=1$. The other elements of $\mathcal{P}$, characterized by $\Tr \brho^2<1$, are called  {\it statistical mixture} or {\it mixed state}. The quantity $\Tr \brho^2$ is thus called the {\it purity} of the state $\brho$.

Pure states are particularly important in quantum mechanics. Indeed, the description of quantum systems in terms of pure density matrices is completely equivalent to the one provided by the state vector (or wavefunction) $|\psi\rangle$. Pure states are, indeed, projectors on one dimensional subspaces of $\mathcal{H}$, thus they uniquely determine a state vector, up to a physically irrelevant global phase  (\cite{breuer-petr,nielsen-chuang}).
In the language of quantum mechanics, we say that a quantum system is {\it closed} if it is isolated from other systems. 
The evolution of an
closed quantum system  is described by the {\it Liouville-von Neumann equation} 
$\dot{\brho} = -i [H, \brho]$,
where $H$ is a Hermitian matrix, called the {\it Hamiltonian} of the system, representing the internal energy of the system. 
 By adding external control fields (such as tunable electromagnetic fields), the perturbed system is governed by a new Hamiltonian which, in most relevant physical situations, can be written as $H(u) = H_0 + \sum_i u_i H_i$, where $H_0$ still represents the internal energy  of the unperturbed system and $H_i$ are associated with the external fields.
In the literature, such controls are usually called {\it coherent controls} (\cite{tutorial-QC}).

With a little abuse of notation, in the field of quantum control we say that a system is closed if its dynamics are described by the Liouville equation, even in presence of interactions with external fields.

If the dependence of $H$ on time is regular enough, $-iH(t)$ is the generator of a unitary evolution group $U_{t,0}$ and, for every $\brho_0\in\mathcal{P}$, the solution of the Liouville equation with initial condition $\brho(0)=\brho_0$ can be written as $\brho(t)=U_{t,0}\brho_0U^{\dagger}_{t,0}$. We remark that the Liouville equation conserves the spectrum of $\brho$; as a consequence, $\Tr(\brho^2(t))$ is constant. In other words, coherent control conserves the purity of a state.

When a quantum system interacts with a surrounding environment,
its evolution is no longer unitary
and reversible, and the general formalism of open quantum
systems is required (\cite{breuer-petr}); under some hypothesis on the environment (such as, Markovianity) the evolution of the quantum system
can be described by the {\it Gorini-Kossakowski-Sudarshan-Lindblad master equation} (see~\cite{gorini-kos-sud,lindblad})
\begin{equation} \label{eq: GKLS}
\dot{\brho} = -i [H, \brho] + \mathcal{L}_{D}(\brho),
\end{equation}
where 
the operator $\mathcal{L}_{D}$ can be written as 
\begin{equation}\label{eq: lind}
\mathcal{L}_{D}(\brho) = \sum_{k=1}^{N^2-1} L_k \brho L^{\dagger}_k - \frac{1}{2} L^{\dagger}_k  L_k  \brho- \frac{1}{2} \brho L^{\dagger}_k L_k,
\end{equation}
the $L_k$ being square matrices called {\it jump (or noise) operators}. In general, the number and the  choice of the jump operators describing the same operator $\mathcal{L}_D$ is not unique; nevertheless, there is always a choice of at most $N^2-1$ operators representing it.

 When $\mathcal{L}\neq 0$, the evolution governed by \eqref{eq: GKLS} preserves the trace and the positivity of the density matrix, but it is no more unitary and isospectral.
 For this reason, usually the term $\mathcal{L}$ is called {\it dissipation}.
 
 Thanks to Choi-Kraus' Theorem (\cite{kraus}), for every $t\geq 0$ the solution of \eqref{eq: GKLS}  with initial condition $\brho_0$ can be written as $\brho(t)=\sum_i M_{i,t}\brho_0M_{i,t}^{\dagger}$, where $\{M_{i,t}\}_i$ is a  family of matrices on $\mathbb{C}^N$ such that  $\sum_i M_{i,t}^{\dagger}M_{i,t}=\id_N$; as for the jump operators, in general, this family is not uniquely determined.


\subsection{Coherence vector representation}

Coherence vector representation is a well known tool in quantum control, that permits to write equation \eqref{eq: GKLS} as a linear differential equation on some real linear space. 
First of all, we endow the space of $N$ dimensional complex square matrices with the Frobenius scalar product $\llangle A,B\rrangle=\Tr\big(A^{\dagger}B\big)$; noticing that every density matrix can be written as $\frac{1}{N}\id_N+\hat{\rho}$, with $\hat{\rho}\in\mathfrak{her}_0(N)$, we choose an orthonormal basis $\{\Lambda_1,\ldots,\Lambda_{N^2-1}\}$ of $\mathfrak{her}_0(N)$ and we  define the map 
$\Phi:\mathcal{P}\to \mathbb{R}^{N^2}$ as 
\[
\Phi(\brho)=\big(\frac{1}{\sqrt{N}},\Tr(\brho \Lambda_1),\ldots,\Tr(\brho \Lambda_{N^2-1})\big).
\]
The map $\Phi$ is called {\it coherence representation} of $\brho$  (\cite{altafiniJMP,dirr-helmke}).
Setting $\overline{\vv^{\rho}}=\Phi(\brho)$ and $\overline{\vv^{\rho}}=(\frac{1}{\sqrt{N}},\vv^{\rho})$, we call $\vv^{\rho}\in\mathbb{R}^{N^2-1}$ the {\it vector of coherences} or {\it Bloch vector} of $\brho$. It is easy to prove that $\Tr(\brho^2)=\|\overline{\vv^{\rho}}\|^2$, which in particular implies that $\|\vv^{\rho}\|^2\leq 1-1/N$. 

 Let us call $K$ the image of $\mathcal{P}$ via the coherence representation; it is easy to see that $\Phi$ is an isomorphism between $\mathcal{P}$ and $K$. 
The properties of $K$ are well-known in the literature; we are resuming them in the following proposition.
 
 \begin{proposition}[\cite{altafiniJMP,avron,dirr-helmke,kimura}]
Let $\boldsymbol{q}_0=(\frac{1}{\sqrt{N}},0,\ldots,0)$, and consider the affine hyperplane $W_0=\{(\frac{1}{\sqrt{N}},\vv) : \vv\in  \mathbb{R}^{N^2-1}\}$.
$K$  
  is a compact convex 
  neighborhood of $\boldsymbol{q}_0$ in $W_0$.
  $\partial K$ is the set of all $\overline{\boldsymbol{v^{\rho}}}$ corresponding to singular density matrices.
\end{proposition}

If $N=2$, then $K$ coincides with the set $\{(\frac{1}{\sqrt{2}},\vv): \|\vv\|\leq \frac{1}{\sqrt{2}}\}$ and can be identified with the three-dimensional ball of radius $\frac{1}{\sqrt{2}}$ (called the {\it Bloch ball}); all vectors in the surface of the Bloch ball  correspond to pure states. 
For $N>2$, the positivity of $\brho$ adds some constraints on $K$, which is strictly contained in the set $\{(\frac{1}{\sqrt{N}},\vv): \|\vv\|\leq \sqrt{1-\frac{1}{N}}\}$; the vectors corresponding to pure states constitute a proper subset of $\partial K$ (more precisely, they correspond to the extreme points of 
$K$).

In coherence representation, equation~\eqref{eq: GKLS} becomes $\dot{\overline{\boldsymbol{v^{\rho}}}}=\big(\Hse{H}+\Hse{\mathcal{L}}\big)\overline{\boldsymbol{v^{\rho}}}$, $\Hse{H}$ and $\Hse{\mathcal{L}}$ being respectively the representations of the operators $-i[H,\cdot]$ and $\mathcal{L}(\cdot)$. It is worth notice (see \cite{altafiniJMP,dirr-helmke} for more details) that, due to the fact that the evolution induced by \eqref{eq: GKLS} is trace-preserving, the matrices  $\Hse{H}$ and $\Hse{\mathcal{L}}$
 always have the block forms 
 \begin{equation} \label{eq: GKLS mat}
 \Hse{H}=
\left(\begin{smallmatrix}
      0 & 0\\
      0 & \Hqu{}
      \end{smallmatrix}
\right)\qquad 
\Hse{\mathcal{L}}=
\left(\begin{smallmatrix}
      0 & 0\\
      \bvz & \widehat{D}
      \end{smallmatrix}
\right),\end{equation} 
where $\Hqu{}\in \mathfrak{so}(N^2-1)$, $\bvz\in \mathbb{R}^{N^2-1}$ and $\widehat{D}$ is a square matrix of dimension $N^2-1$. 

The coherence representation of equation~\eqref{eq: GKLS} subject to coherent control is 
\begin{equation} \label{eq: coherent controlled GKLS}
\dot{\overline{\boldsymbol{v^{\rho}}}}=\big(\Hse{H_0}+\sum_j u_j \Hse{H_j} +\Hse{\mathcal{L}}\big)\overline{\boldsymbol{v^{\rho}}}    .
\end{equation}

The controllability properties of equation~\eqref{eq: coherent controlled GKLS} have been studied in  \cite{altafiniJMP}, where, in particular, it is stated that the system is never STLC. We recall moreover the following result.
\begin{theorem} \label{th: boundary unreachable}
Let $\bq_0\in \mathring{K}$ and $\bq_1 \in  \partial K$. There is no essentially bounded control function 
that steers $\bq_0$ to $\bq_1$
in finite (positive) time.
\end{theorem}


\subsection{The system} \label{sec: system}
In this paper, we focus on coherently controlled composite open quantum systems,  composed by two interacting qubits (called $A$ and $B$); we recall that a qubit is a quantum system living in a two-dimensional Hilbert space.
According to quantum mechanics (\cite{breuer-petr,haroche,nielsen-chuang}), the total system $A+B$ evolves on the tensor product of two Hilbert spaces, 
$\mathcal{H} = \mathcal{H}_A \otimes  \mathcal{H}_B$, with $\mathcal{H}_A=\mathcal{H}_B=\mathbb{C}^2$.
Also, we recall that the state of each subsystem (denoted respectively with $\rho_A$ or $\rho_B$) can be
extracted from the state $\brho$ of the total system, by means of the  partial trace:  indeed, 
we stress that, if $\brho$ is a density matrix (i.e. Hermitian positive semi-definite and of trace one) on $\mathcal{H}$, then $\rho_A=\Tr_B\brho$ and $\rho_B=\Tr_A\brho$ are density matrices  on $\mathcal{H}_A$ and $\mathcal{H}_B$, respectively. 
Physically, taking the partial trace on $A$  can be interpreted as ``averaging'' on the information about $A$, so that the reduced state $\rho_B$ describes the state of the subsystem $B$. Analogous considerations hold for $\rho_A$.

The system $A+B$ evolves according to equation~\eqref{eq: GKLS}.  We remark that we can always assume that the matrices $H$ and $L_k$ are traceless, as adding them multiples of the identity leaves invariant the right-hand side of \eqref{eq: GKLS}. 
Then, 
$H$ 
can be uniquely written as $H= H_{A}+H_{I}+H_{B}$, where $H_{A}=\Hdue{A}\otimes \id_2$,  $H_{B}=\id_{n_A}\otimes \Hdue{B}$, with $\Hdue{A},\Hdue{B}\in \mathfrak{her}_0(2)$ and
$H_I\in \mathfrak{her}_0(2)\otimes \mathfrak{her}_0(2)$.
Actually, $H_I$ denotes the interaction between the subsystems $A$ and $B$, while $H_A$ and $H_B$  the unitary free evolution of the two subsystems.

We make these further assumptions on the dynamics:
\begin{itemize}
\renewcommand{\labelitemi}{$\cdot$}
 \item the subsystem $B$ does not interact directly with the  environment; in particular, this implies that the jump operators can be taken of the form $L_k=\ell_k\otimes\id_2$, where the $\ell_k$'s are traceless matrices on $\mathcal{H}_A$.
 \item the control directly affects only the subsystem $A$; in particular, the Hamiltonian $H$ can be written as 
\[
H=\big(\Hdue{A_0}+\sum_{j=1}^3u_j\Hdue{A_j}\big)\otimes \id_2+H_I+\id_2\otimes \Hdue{B},
\]
where the matrices $\Hdue{A_j},\ j\geq 0$, $H_I$ and $ \Hdue{B}$ are constant and $u_j$ are functions of time.
\end{itemize}

It is clear that, if $H_I=0$, then the two subsystems are completely independent, and, in particular, the evolution of $B$ is not influenced by the control. 

As already said, in \cite{altafiniJMP} it is proved that the coherently controlled Lindblad equation is never STLC and Theorem~\ref{th: boundary unreachable} states that some transitions cannot be realized in finite time.
Nevertheless, it is still interesting to understand what can be said about the controllability of the state of the subsystem $B$ only.
 In particular, the main questions that one may ask include
 \begin{enumerate}
 \item \label{q: constant purity} is it possible to ``protect'' the subsystem $B$ from dissipation, that is, to implement a unitary dynamics on $B$, at least on some submanifolds (for instance, the submanifold  $\Tr \rho_B^2=1$)? 
  \item \label{q: control} is it possible to ``control'' the state $\rho_B$, regardless of $\brho$?
 \end{enumerate}

  Question \ref{q: control}) needs to be further clarified, as several notions of ``partial controllability'' are possible (see for instance \cite{alb-rom-dal}); in any of its declination, tackling the issue is a very hard task, even in the simplest case of two qubits. 
  
  Question \ref{q: constant purity}) seems to be more affordable. In this paper, we provide a first step towards the answer: for three particular choices of the interaction $H_I$ (the well-known dispersive and resonant couplings, very common in experimental set-ups), and a generic choice of the dissipative term, we show that the trajectories keeping $\rho_B$ pure are trivial, or follow a free evolution which is not affected by the controls. In our opinion, more general interaction would lead to similar results.

\section{On the structure of $K$} \label{sec: K}

In the case of two interacting qubits, we choose the following orthonormal basis of $\mathfrak{her}_0(4)$
\begin{align}
\Lambda_i&=\frac{1}{2}\sigma_i\otimes \id_{2},\qquad i=1,2,3\\ 
\Lambda_{3i+j}&=\frac{1}{2}\sigma_i\otimes \sigma_j,\qquad i=1,2,3, \ j=1,2,3,\\ 
\Lambda_{12 +j}&=\frac{1}{2}\id_2\otimes \sigma_j,\qquad j=1,2,3,
\end{align}
where $\sigma_i$, $i=1,2,3$, denote the {\it Pauli matrices}
\begin{equation}\label{eq: Pauli}
\sigma_1=\begin{pmatrix}
          0 & 1 \\ 1 & 0
         \end{pmatrix},\
\sigma_2=\begin{pmatrix}
          0 & -i \\ i & 0
         \end{pmatrix}, \
         \sigma_3=\begin{pmatrix}
          1 & 0 \\ 0 & -1
         \end{pmatrix}. 
\end{equation}
For such system, in the following we will adopt also the more intuitive notation
\begin{equation}\label{eq: bloch vec 2x2}
\boldsymbol{v^{\rho}}=\big(\frac{1}{2},\boldsymbol{v^A},\boldsymbol{v^{AB}},\boldsymbol{v^B}\big),
\end{equation}
with $\boldsymbol{v^A},\boldsymbol{v^B}\in \mathbb{R}^3$, and $\boldsymbol{v^{AB}}\in \mathbb{R}^{9}$.

Thanks to the peculiarity of the Pauli matrices (in particular, the fact that $\sigma_i\sigma_j+\sigma_j\sigma_i=2\delta_{ij}\id_2$), the matrices $\Hqu{}$ and $\widehat{D}$ and the vector $\bvz$ in
equation \eqref{eq: GKLS mat} have the following block structure
\begin{align} \label{eq: Ham coh}
 \Hqu{}&=\left(  
 \begin{array}{c|c|c}
  \Htre{A} & \widehat{H}_{It} & {\bf 0}\\
  \hline
  \widehat{H}_{Il}& \Htre{A}\otimes \id_{3} + \id_{3}\otimes \Htre{B} & \widehat{H}_{Ir}\\
  \hline
  {\bf 0}& \widehat{H}_{Ib} &  \Htre{B} \\ 
 \end{array}
 \right)\\
\widehat{D}&=\left(
\begin{array}{c|c|c}
 \hat{d} &{\bf 0}& {\bf 0}\\ \hline
  {\bf 0} & \hat{d}\otimes \id_3& \bvz\otimes \id_3\\ \hline
 {\bf 0} & {\bf 0}& {\bf 0}
\end{array}
\right)\\
\bvz&=({\bvz}_1,{\bvz}_2,{\bvz}_3,0,\ldots,0)
\end{align}
(more details on the objects appearing here above are provided in Appendix\ref{app: computations}).

In order to answer to question~\ref{q: constant purity}, 
 we look for a characterization of states $\brho$ such that $\Tr_A \brho$ is pure. 
 \begin{proposition}\label{prop: pp bordo} \label{prop: fact}
 Consider $\brho\in \mathcal{P}$ such that $\Tr\rho_B^2=1$, 
 where $\rho_B=\Tr_A\brho$. Then $\Phi(\brho)\in \partial K$
%
%
%
%
and there exist $\boldsymbol{v}^A,\boldsymbol{v}^B\in \mathbb{R}^3$
with $\|\vv^A\|^2\leq \frac{1}{4}$ and $\|\vv^B\|^2= \frac{1}{4}$
such that
\begin{equation} \label{eq: fact}\Phi(\brho)=
 \big(\frac{1}{2},\vv^A,2 \vv^A\otimes \vv^B,\vv^B\big).
 \end{equation}
 \end{proposition}

 \begin{proof}
 The first claim is an easy consequence of \cite[Proposition~2.1]{Lin2016}. 
 
Let us now prove the second part. 
The bound on $\|\vv^A\|$ and the value of $\|\vv^B\|$ yield from the fact that $\Tr(\Tr_B\brho)^2=\frac{1}{2}+2\sum_{j=1}^{3}(v^A_j)^2$ and
$\Tr(\Tr_A\brho)^2=\frac{1}{2}+2\sum_{j=1}^{3}(v^B_j)^2$.

Moreover 
\begin{align} \label{eq: derivata traccia}
&\frac{1}{4} \frac{d}{dt}\Tr\rho_B^2=
\langle \boldsymbol{v}^{B}, \widehat{h}_{B}\boldsymbol{v}^{B}\rangle+\langle \boldsymbol{v}^{B}, \widehat{H}_{Ib}\boldsymbol{v}^{AB}\rangle\\
&= \sum_{r,s=1}^{3} \lambda_{rs} \langle  \label{eq: partial trace variation}
\boldsymbol{v}^{B}, T_s (v^{AB}_{(r-1)+1},v^{AB}_{(r-1)+2},v^{AB}_{(r-1)+3}) \rangle,
\end{align}
where we used equation~\eqref{eq: block H_ib} and the fact that $\widehat{h}_{B}$ is antisymmetric. 

By hypothesis, $\Tr \rho_B^2$ cannot increase with time. 
On the other hand, as $\{T_1,T_2,T_3\}$ constitute a basis of $\mathfrak{so}(3)$, and equation~\eqref{eq: derivata traccia} holds for any $H_I$,
 $\langle 
\boldsymbol{v}^{B}, A (v^{AB}_{(r-1)+1},v^{AB}_{(r-1)+2},v^{AB}_{(r-1)+3}) \rangle$
 must be zero for every antisymmetric matrix $A$ and for $r=1,2,3$.
Thus, the vectors $(v^{AB}_{(r-1)+1},v^{AB}_{(r-1)+2},v^{AB}_{(r-1)+3})$ must be collinear to $\boldsymbol{v}^{B}$, and
$\brho$ can be written as
\begin{align}
\brho&=\frac{1}{4}\id_4+\sigma_A\otimes\frac{\id_2}{\sqrt{2}}\\
&+\frac{1}{2}\big(\gamma \sigma_1+\beta \sigma_2+\theta\sigma_3+\id_2\big)\\ 
&\qquad \qquad \otimes \big(\boldsymbol{v_1^B} \sigma_1+\boldsymbol{v_2^B} \sigma_2+\boldsymbol{v_3^B}\sigma_3),
\end{align}
for some $\sigma_A=\frac{1}{\sqrt{2}}\big(\boldsymbol{v^A_1}\sigma_1+\boldsymbol{v^A_2}\sigma_2+\boldsymbol{v^A_3}\sigma_3\big)$.
Let $\{\psi_1,\psi_2\}$ be the (normalized) eigenvectors of  $\boldsymbol{v_1^B} \sigma_1+\boldsymbol{v_2^B} \sigma_2+\boldsymbol{v_3^B}\sigma_3$, relative to the eigenvalues $\mu_1=-\frac{1}{2}$ and $\mu_2=\frac{1}{2}$, respectively, and let  $\{\phi_1,\phi_2\}$ be any (possibly different) orthonormal basis  of $\mathbb{C}^2$; 
let $\lambda$ be an eigenvalue of $\brho$, and $\varphi$ a corresponding eigenvector, which can be written as  
$\varphi=\sum_{i,j=1}^2 a_{ij} \phi_i\otimes \psi_j$.
Set $A_j=\frac{\sigma_A}{\sqrt{2}}+\frac{\mu_j}{2}(\gamma \sigma_1+\beta \sigma_2+\theta\sigma_3+\id_2)+\frac{1}{4}\id_2$, for $j=1,2$. We remark that the spectrum of $\brho$ is given by the union of the spectra of $A_1$ and $A_2$.
In particular,  the eigenvalues of $A_1$ are given by  
\[
\nu_{\pm}= \pm\sqrt{\big(\frac{\boldsymbol{v^A_1}}{2}-\frac{\gamma}{4}\big)^2+\big(\frac{\boldsymbol{v_2^A}}{2}-\frac{\beta}{4}\big)^2+\big(\frac{\boldsymbol{v^A_3}}{2}-\frac{\theta}{4}\big)^2}.
\]
As they belong also to the spectrum of $\brho$, which is positive semidefinite, then it must be $\boldsymbol{v^A_1}=\frac{\gamma}{2}$, $\boldsymbol{v^A_2}=\frac{\beta}{2}$, and $\boldsymbol{v^A_3}=\frac{\theta}{2}$, and the proposition is proved.
\end{proof}


 \begin{remark}
Proposition~\ref{prop: fact} can be easily generalized to the case in which the subsystem $A$ has (complex) dimension $n_A\geq2$.
\end{remark}

\section{First answers to question \ref{q: constant purity})} \label{sec: results}

Proposition~\ref{prop: fact} imposes a constraint on the structure of states $\brho$ whose reduction to $B$ is pure, and can be thus exploited to study the controllability of the reduced states. 
First of all, together with Theorem~\ref{th: boundary unreachable}, it yields the following fact. 
 \begin{corollary} \label{cor: ppure infinite time}
 Let $\brho_0,\brho_T\in \mathcal{P}$ such that $ \Phi(\brho_0)$ belongs to the interior of $K$ and $\Tr(\Tr_A\brho_T)^2=1$. Then, there is no essentially bounded control function that can send $\brho_0$ to $\brho_T$ in finite time.   
 \end{corollary}

\medskip 
In other to give partial answers to Question~\ref{q: constant purity}, we study the following problem.
\begin{prob}
Let $\brho_0$ such that $\Tr\big(\Tr_A \brho_0)^2=1$, and call $\brho(t)$ the solution at time $t$ of~\eqref{eq: GKLS}, under the assumption in Section~\ref{sec: system}, with $\brho(0)=\brho_0$.  Is it possible to find $\epsilon>0$ and a control function defined on $[0,\epsilon]$  such that $\Tr\big(\Tr_A \brho(t))^2=1$  for $t\in [0,\epsilon]$? 
\end{prob}

Assume that it is true, for some  piecewise-$C^{\infty}$ control function $\hat{\bu}:[0,\epsilon]\to\mathbb{R}^3$, 
Set $\Phi(\brho(t))=\big(\frac{1}{2},\boldsymbol{v^A}(t),\boldsymbol{v^{AB}}(t),\boldsymbol{v^B}(t)\big)$. By Proposition~\ref{prop: fact}, if $\Tr \rho^2_B(t)=1$ for $t\in[0,\epsilon]$, then $\boldsymbol{v^{AB}}(t)=2\boldsymbol{v^A}(t)\otimes\boldsymbol{v^{B}(t)}$. In particular, for every $k\geq 1$, we have
\begin{equation} \label{eq: derivatives}
\frac{d^k}{dt^k}\boldsymbol{v^{AB}}(t)=2\sum_{j=0}^k \begin{pmatrix} 
                                                   k \\ j
                                                  \end{pmatrix}
\frac{d^j\boldsymbol{v^{A}}}{dt^j}\otimes \frac{d^{k-j}\boldsymbol{v^{B}}}{dt^{k-j}}.\end{equation}
 
 In the following, we will use equations~\eqref{eq: derivatives} to find out the controls satisfying the claim (if any), for different expressions of the interaction $H_I$. In order to do it, we 
 define the vector
 \[\boldsymbol{w}(t)=\boldsymbol{\dot{v}^{AB}}(t)-2\sum_{j=0,1} \left(\begin{smallmatrix} 
                                                   k \\ j
                                                  \end{smallmatrix}\right)
\frac{d^j\boldsymbol{v^{A}}}{dt^j}\otimes \frac{d^{k-j}\boldsymbol{v^{B}}}{dt^{k-j}}.\]

  \medskip
 Without loss of generality, we can choose bases on $\mathcal{H}_A$ and $\mathcal{H}_B$ such that  
\[
H_{A_0}=\omega_a \sigma_3\quad
H_{B_0}=\omega_b \sigma_3\quad
\]
Eventually performing a linear transformation in the control space, we also assume that $H_{A_i}= \sigma_i$ for $i=1,2,3.$

\smallskip
\noindent
{\bf ``Dispersive'' coupling: $H_I=g\sigma_3\otimes\sigma_3$.}
Computing the vector $\boldsymbol{w}$, 
we notice that 
\begin{gather}
w_3=g\boldsymbol{v_2^A}(1-4(\boldsymbol{v_3^B})^2)\\ w_8=g\boldsymbol{v_2^B}(1-4(\boldsymbol{v_3^A})^2),
\end{gather}so that equation~\eqref{eq: derivatives} is satisfied only
if $(\boldsymbol{v_3^B})^2$ is identically equal to $\frac{1}{4}$ (which implies $\boldsymbol{v_1^B}\equiv \boldsymbol{v_2^B}\equiv 0$) 
or if $(\boldsymbol{v_3^A})^2\equiv \frac{1}{4}$ (which implies $\boldsymbol{v_1^A}\equiv \boldsymbol{v_2^A}\equiv 0$).

We remark that the first scenario corresponds to freezing $\rho_B$ to the state $\frac{1}{2}\big(\id_2+\sigma_3\big)$ or to the state $\frac{1}{2}\big(\id_2-\sigma_3\big)$, that is, $\rho_B$ is constant. In particular, by computation we notice that, for this choice of $H_I$, equation~\eqref{eq: coherent controlled GKLS} has the block-triangular form  
\[
\begin{pmatrix}
\dot{z}_1 \\ \dot{z}_2 \\ \dot{v}_3^B
\end{pmatrix}=\begin{pmatrix}
C_{11} & C_{12} & C_{13}\\
0 & C_{22} & 0\\
0 & 0 & 0
\end{pmatrix}
\begin{pmatrix}
z_1\\z_2\\ \boldsymbol{v_3^B}
\end{pmatrix}
\]
with $z_2=(\boldsymbol{v^{AB}_1},\boldsymbol{ v^{AB}_2 },\boldsymbol{ v^{AB}_4},\boldsymbol{ v^{AB}_5},\boldsymbol{ v^{AB}_7},\boldsymbol{ v^{AB}_8},\boldsymbol{ v^{B}_1},\boldsymbol{ v^{B}_2})$.
If $\boldsymbol{v^{B}_3}=|\frac{1}{2}|$ and $\rho_B$ is pure at $t=0$, then $z_2(0)=0$, so that $z_2(t)$ is zero for all $t\geq 0$, and every value of the control. In particular, 
the submanifolds of $\mathcal{P}$ 
\begin{gather}
\{\rho_A\otimes\frac{1}{2}\big(\id_2+\sigma_3\big): \rho_A=\rho_A^{\dag},\rho_A\geq0,\Tr \rho_A=1\}\\ \{\rho_A\otimes\frac{1}{2}\big(\id_2-\sigma_3\big): \rho_A=\rho_A^{\dag},\rho_A\geq0,\Tr \rho_A=1\},   
\end{gather}  
are invariant for equation~\eqref{eq: GKLS}, for every choice of the control and for any dissipation.

Let us now consider the second case; first of all, by \eqref{eq: coherent controlled GKLS}, we notice that $\boldsymbol{\boldsymbol{v^A_3}}$ can be constantly equal to $|\frac{1}{2}|$ only if  $\frac{{\bvz}_3}{2}+\boldsymbol{v_3^A} \hat{d}_{33}=0$, which restrict the class of dissipative terms $\mathcal{L}$ that allow such behavior.
Some examples of such operators are the so-called {\it amplitude damping channels} (\cite{nielsen-chuang}), that is,
associated respectively to the jump operators $\sigma_+\otimes \id_2$ or $\sigma_-\otimes \id_2$, where $\sigma_\pm=\sigma_1\pm i\sigma_2$.

On the other hand, taking into account the fact that the state is factorized, the equation for $\boldsymbol{v}^B$ becomes 
\[
\begin{pmatrix}
\boldsymbol{\dot{v}_1^B}\\ \boldsymbol{\dot{v}_2^B}\\ \boldsymbol{\dot{v}_3^B}
\end{pmatrix}\!\!=\!
\begin{pmatrix}
0  & \!\!-\omega_b-2g \boldsymbol{v_3^A} & 0\\
\omega_b +2g \boldsymbol{v_3^A}& 0 & 0\\
0 & 0 &0
\end{pmatrix}\!\!
\begin{pmatrix}
\boldsymbol{v_1^B}\\ \boldsymbol{v_2^B}\\\boldsymbol{v_3^B}
\end{pmatrix}\!\!,
\]
which is unaffected by the control.

\smallskip
\noindent
{\bf ``Resonant'' coupling: $H_I=g(\sigma_+\otimes\sigma_-+\sigma_-\otimes \sigma_+)$.} 
First of all, we remark that we can write $H_I=\frac{g}{2}(\sigma_1\otimes\sigma_1+\sigma_2\otimes \sigma_2)$.

As we did above, we try to find conditions that guarantee that  $\boldsymbol{w}$ is null.
As
\[\Scale[0.9]
{\boldsymbol{w}=g
\begin{pmatrix}
-4 (v^A_3v^B_1 v^B_2 + 
   v^A_1 v^A_2 v^B_3) 
   \\ -(v^A_3 (-1 + 4 (v^B_2)^2) + v^B_3 (1- 4 (v^A_1)^2) 
   \\ (4 v^A_1 (v^A_2v^B_1-v^A_1 v^B_2) + 
   v^B_2 (1 - 4 v^A_3  v^B_3) )
   \\ (v^A_3 (-1 + 4 (v^B_1)^2) + v^B_3(1 - 
   4 (v^A_2)^2)) 
   \\ 4 (v^A_3v^B_1 v^B_2 + v^A_1 v^A_2 v^B_3) 
   \\ ( 4 v^A_2(v^A_2v^B_1-v^A_1 v^B_2) + 
  v^B_1 (-1  + 4 v^A_3 v^B_3)) 
   \\ ( 4v^B_1 (v^A_2v^B_1-v^A_1 v^B_2)  + 
    v^A_2(1-4 v^A_3 v^B_3)) 
   \\ (- 4 v^B_2 (v^A_2v^B_1-v^A_1 v^B_2) + 
   v^A_1 (-1  + 4 v^A_3 v^B_3)) 
   \\ 4 (v^A_2v^B_1 - v^A_1 v^B_2) (v^A_3 - v^B_3),
\end{pmatrix}}\]
we deduce that, if $\boldsymbol{w}$ is null, then  $\boldsymbol{v^A_2 v^B_1} = \boldsymbol{v^A_1 v^B_2}$ and/or $\boldsymbol{v^A_3}=\boldsymbol{v^B_3}$.

In the first case, plugging the equality into the expression of $\boldsymbol{w}$ and setting it to zero, we obtain
\[
\begin{cases}
\boldsymbol{v^A_1} (1-4 \boldsymbol{v^A_3} \boldsymbol{v^B_3})=0\\
\boldsymbol{v^A_2} (1-4 \boldsymbol{v^A_3} \boldsymbol{v^B_3})=0\\
\boldsymbol{v^B_1} (1-4 \boldsymbol{v^A_3} \boldsymbol{v^B_3})=0\\
\boldsymbol{v^B_2} (1-4 \boldsymbol{v^A_3} \boldsymbol{v^B_3})=0.
\end{cases}
\]
 These four equations are satisfied if $\boldsymbol{v^A_1}=\boldsymbol{v^A_2}= \boldsymbol{v^B_1} =\boldsymbol{v^B_2}=0$, or if $\boldsymbol{v^A_3} \boldsymbol{v^B_3}=1/4$; in both cases, due to the fact that $\rho_B$ is pure and to the constraints on the length of $\vv^A$ and $\vv^B$, we obtain that $|\boldsymbol{v^A_3}|=|\boldsymbol{v^B_3}|=1/2$ and that also $\rho_A$ is a pure state (and, as a consequence of the ``factorized structure'', the whole state $\brho$ is pure). 
 
 If instead $\boldsymbol{v^A_3}=\boldsymbol{v^B_3}$, setting to zero the second and the fourth  components of $\boldsymbol{w}$, we obtain that $|\boldsymbol{v^A_1}|=|\boldsymbol{v^B_2}|$ and $|\boldsymbol{v^A_2}|=|\boldsymbol{v^B_1}|$, which again implies that $\rho_A$ is a pure state.
 
 Summing up, in presence of resonant coupling, it is not possible to keep the partial state $\rho_B$ pure if the whole state $\brho$ itself is not kept pure by the evolution.  
 On the other hand, as
 $\frac{d}{dt}\Tr\brho^2=2\Tr(\brho\mathcal{L}(\brho))$, it follows that the state $\rho_B$ can be kept pure only if $\brho$ evolves in the set $\Tr(\brho\mathcal{L}(\brho))=0$; for ``factorized states'' (that is, of the form \eqref{eq: fact}), this happens only for states such that $\frac{\mathbf{v}_0}{2}+\hat{d} \boldsymbol{v^A}=0$; depending on the particular choice of $\mathcal{L}$, this equation may not have solutions $\boldsymbol{v^A}$ of norm $1/2$. \\

 \noindent
 $\boldsymbol{H_I=g\sigma_3\otimes\sigma_1}$. We finally discuss a further case in which the interaction does not commute with the free Hamiltonian $H_B$.
 
 Setting $\boldsymbol{w}\equiv 0$, we find the following constraints: either $\boldsymbol{v_2^B}\equiv \boldsymbol{v_3^B} \equiv 0$ and $(\boldsymbol{v^B_1})^2\equiv 1/4$, or $\boldsymbol{v^A_1}\equiv \boldsymbol{v^A_2}\equiv 0$ and $(\boldsymbol{v^A_3})^2\equiv 1/4$. 
 
 The first case does not correspond to an admissible solution of the control system: indeed, by computations, it is possible to see that no control can keep $\brho$ in a state of the form $\rho_A\otimes \frac{1}{2}\big(\id_2\pm\sigma_1\big)$ on a nonzero time interval. 
 
Let us now look for admissible trajectories along which $\boldsymbol{v^A_1}\equiv \boldsymbol{v^A_2}\equiv 0$ and $(\boldsymbol{v^A_3})^2\equiv 1/4$.
First of all, as 
$\boldsymbol{\dot{v}^A_3}=\frac{{\bvz}_3}{2}+\boldsymbol{v^A_3} \hat{d}_{33}+\boldsymbol{v^A_1}(\hat{d}_{31}-u_2)+\boldsymbol{v^A_2}(\hat{d}_{32}+u_1)$, such trajectories are admissible only if $\boldsymbol{v_3^A}\hat{d}_{33}+\frac{{\bvz}_3}{2}=0$, that is true for some particular dissipation terms only (as we already saw, the amplitude damping  channels satisfy such a constraint). 

Inspecting the differential equations for $\boldsymbol{v^A_1}$ and $\boldsymbol{v^A_2}$, we see that they stay constant only for the choice of the control
\[
u_1 =\frac{{\bvz}_2+2 \boldsymbol{v^A_3} \hat{d}_{23}}{2 \boldsymbol{v^A_3}} \qquad  
u_2=-\frac{{\bvz}_1+2 \boldsymbol{v^A_3} \hat{d}_{13}}{2 \boldsymbol{v^A_3}}.
\]
Moreover, as this guarantees $\boldsymbol{w}\equiv 0$, that is, $\boldsymbol{v}_{\rho}$ has the form \eqref{eq: fact}, we can substitute the values of $\boldsymbol{v}^{AB}$ into the differential equations for $\boldsymbol{v}^B$, getting
\[
\begin{pmatrix}
\boldsymbol{\dot{v}_1^B}\\ \boldsymbol{\dot{v}_2^B}\\ \boldsymbol{\dot{v}_3^B}
\end{pmatrix}=
\begin{pmatrix}
0  & -\omega_b & 0\\
\omega_b & 0 & -2g \boldsymbol{v_3^A}\\
0 & 2g \boldsymbol{v_3^A} &0
\end{pmatrix}
\begin{pmatrix}
\boldsymbol{v_2^B}\\ \boldsymbol{v_2^B}\\\boldsymbol{v_3^B}
\end{pmatrix},
\]
that is, the dynamics of $\rho_B$ are protected from dissipation, but the control does not affect them.

 \section{Conclusions}
In this paper, we presented a first analysis on the indirect controllability properties of a 2-qubit system, in the case in which the ancilla is subject to dissipation.

 First of all, we observed that states $\brho$ such that their reduction $\Tr_A \brho$ is pure are not reachable (in finite time) from the interior of the space $\mathcal{P}$, thus obtaining a first obstruction to the indirect controllability of the system; it would be interesting to investigate if such states are reachable from any other point of the boundary of $\mathcal{P}$. 

 We then focus on the possibility of preserve the subsystem $B$  from dissipation, that is, to find admissible trajectories $\brho(t)$ such that $\Tr_A \brho(t)$ is pure. 
We investigated three particular cases of interaction (among them, the well known dispersive and resonant couplings), and we found that the only admissible trajectories are either trivial (i.e. their reduction $\Tr_A \brho(t)$ is constant) or are unaffected by the action of the control. In our opinion,  similar results hold also for other interaction Hamiltonians $H_I$.

 \appendix
 \section{Useful formulas for the two qubits case} \label{app: computations}

First of all we recall that the matrices
\[
T_1=\left(\begin{smallmatrix}
      0 & 0 & 0\\
            0 & 0& -1\\
                  0 & 1 & 0
     \end{smallmatrix}\right)\quad
T_2=\left(\begin{smallmatrix}
      0 & 0 & 1\\
            0 & 0& 0\\
                  -1 & 0 & 0
     \end{smallmatrix}\right)\quad 
     T_3=\left(\begin{smallmatrix}
      0 & -1 & 0\\
            1 & 0& 0\\
                  0 & 0 & 0
     \end{smallmatrix}\right)
\]
are the representations on $\mathbb{R}^3$ of the operators $-i[\tfrac{\sigma_j}{2},\cdot]$, $j=1,2,3$. Then, if $\Hdue{A}=\sum_{j=1}^3 \frac{\alpha_j}{2}\sigma_j$ and  $\Hdue{B}=\sum_{j=1}^3 \frac{\beta_j}{2}\sigma_j$ (we recall that we assumed them traceless), the matrices $\Htre{A}$ and $\Htre{B}$ in \eqref{eq: Ham coh} are simply
$\Htre{A}=\sum_{j=1}^3 \alpha_j T_j$ and $\Htre{B}=\sum_{j=1}^3 \beta_j T_j$.
 
Writing $
    H_I=\frac{1}{2}\sum_{i,j=1}^3 \lambda_{ij} \sigma_i\otimes \sigma_j$, long but easy computations give
 \begin{align}
   \widehat{H}_{It}&=\sum_{j=1}^3 T_j\otimes (\lambda_{j1},\lambda_{j2},\lambda_{j3})\\ 
  \widehat{H}_{Ib}&=\sum_{j=1}^3  (\lambda_{1j},\lambda_{2j},\lambda_{3j})\otimes T_j   \label{eq: block H_ib} 
 \end{align}
and the blocks $\widehat{H}_{Il}$ and $\widehat{H}_{Ir}$ can be recovered by antisymmetry.
 
%
%
%
%
%
%
%
%

\bibliography{kateryna-paper.bbl}{}
\bibliographystyle{plain}

\end{document}